\documentclass[]{interact}

\usepackage{epstopdf}
\usepackage[caption=false]{subfig}

\usepackage{color}
\usepackage[nice]{nicefrac}
\usepackage{floatrow}

\usepackage[numbers,sort&compress]{natbib}
\bibpunct[, ]{[}{]}{,}{n}{,}{,}
\makeatletter
\def\NAT@def@citea{\def\@citea{\NAT@separator}}
\makeatother

\theoremstyle{plain}

\theoremstyle{definition}

\theoremstyle{remark}

\begin{document}


\title{The scaling properties and the multiple derivative of Legendre polynomials}

\author{
\name{G.~M. Laurent and G.~R. Harrison \thanks{CONTACT G.~M. Laurent. Email: glaurent@auburn.edu}}
\affil{Department of Physics, Auburn University, 206 Allison Laboratory, Auburn, AL, 36849, USA}
}

\maketitle

\begin{abstract}
In this paper, we study the scaling properties of Legendre polynomials $P_{n}(x)$. We show that $P_{n}(\lambda x)$, where $\lambda$ is a constant, can be expanded as a sum of either Legendre polynomials $P_{n}(x)$ or their multiple derivatives $d^{k}P_{n}(x)/dx^{k}$, and we derive a general expression for the expansion coefficients. In addition, we demonstrate that the multiple derivative $d^{k}P_{n}(x)/dx^{k}$ can also be expressed as a sum of Legendre polynomials and we obtain a recurrence relation for the coefficients. 
\end{abstract}


\section{Introduction}

The central importance of Legendre polynomials in many fields of pure and applied sciences is undoubtedly well-established. Consequently, their properties have been extensively investigated over several decades \cite{Gradshteyn14, Olver10}, and still remain a matter for numerous studies \cite{Bosch00, Anli07, Antonov10, Szmytkowski11, Dattoli11, Bos17}. To our knowledge, though, the scaling properties of Legendre polynomials, hereafter denoted as $P_{n}(x)$, have not been reported in the literature so far. The purpose of this work is to derive explicit expressions for $P_{n}(\lambda x)$, where $\lambda$ is a constant. Such expressions appear to be very useful for deriving an analytical representation for the projection of spherical harmonics onto a plan. In this paper, we show that $P_{n}(\lambda x)$ can be expanded either as a sum of Legendre polynomials $P_{n}(x)$ or their multiple derivatives:
\begin{equation}
P_{n}(\lambda x)= \sum_{k=0}^{\lfloor n/2 \rfloor} a_{\lambda, n, k} \dfrac{d^{k}}{dx^{k}}P_{n-k}(x)= \sum_{k=0}^{\lfloor n/2 \rfloor} b_{\lambda, n, k} P_{n-2k}(x),
\label{Eq1a} 
\end{equation}
and we derive an expression for the expansion coefficients $a_{\lambda, n, k}$ and $b_{\lambda, n, k}$. In addition, we demonstrate that the multiple ($k^{th}$) derivative of Legendre polynomials of degree $n$ can also be expressed as a sum of Legendre polynomials:
\begin{equation}
\dfrac{d^{k}}{dx^{k}}P_{n}(x)=\sum_{i=0}^{\lfloor (n-k)/2 \rfloor} \alpha_{n-k-2i}P_{n-k-2i}(x),
\label{Eq1b}
\end{equation}
 and we derive a recurrence relation for the coefficients $\alpha_{n-k-2i}$. 

\section{Proof}

\subsection{Scaling properties}

Using Rodrigues' formula, Legendre polynomials $P_{n}(\lambda x)$ of degree $n$ can be written as \cite{Gradshteyn14}:
\begin{equation}
P_{n}(\lambda x)= \dfrac{1}{2^{n}n!} \dfrac{d^{n}}{d(\lambda x)^{n}}\left[(\lambda^{2} x^{2}-1)^{n}\right].
\label{Eq2}
\end{equation}
By recognizing that $d(\lambda x)^{n}=\lambda ^{n}dx^{n}$ and $\lambda^{2} x^{2}-1=\lambda ^{2}(x^{2}-1)+\lambda ^{2}-1$, equation (\ref{Eq2}) can be  rewritten as:
\begin{equation}
P_{n}(\lambda x) = \dfrac{1}{2^{n}n!} \lambda ^{n}\dfrac{d^{n}}{dx^{n}}\left[\left((x^{2}-1)+\dfrac{\lambda ^{2}-1}{\lambda ^{2}}\right)^{n}\right].
\label{Eq3}
\end{equation}
In the previous equation, the polynomial of degree $n$ being derived can be expanded as a sum of $n+1$ terms:
\begin{equation}
\left((x^{2}-1)+\dfrac{\lambda ^{2}-1}{\lambda ^{2}}\right)^{n} = \sum_{k=0}^{n} \alpha_{k}(x^{2}-1)^{k},
\label{Eq4}
\end{equation}
whose expansion coefficients, $\alpha_{k}$, can be related to the product of its unique root $(\lambda^{2}-1)/\lambda^{2}$ using Vieta's formula \cite{Polyanin07}. It can be shown that:
\begin{equation}
\alpha_{k}={n \choose k} \left(\dfrac{\lambda ^{2}-1}{\lambda ^{2}}\right)^{n-k}.
\label{Eq5}
\end{equation}
From equations (\ref{Eq3}), (\ref{Eq4}), and (\ref{Eq5}), it follows that:
\begin{equation}
P_{n}(\lambda x) = \dfrac{1}{2^{n}n!} \dfrac{d^{n}}{dx^{n}} \sum_{k=0}^{n} {n \choose k} \dfrac{(\lambda^{2}-1)^{n-k}}{\lambda ^{n-2k}} (x^{2}-1)^{k}.
\label{Eq6}
\end{equation}
By interchanging the order of the derivation and the sum in equation (\ref{Eq6}) and recognizing that the $n^{th}$ derivative of the polynomial $(x^{2}-1)^{k}$ vanishes for $k$ lower than  $\lceil n/2\rceil$, we obtain:
\begin{equation}
P_{n}(\lambda x) = \sum_{k=\lceil n/2 \rceil}^{n} \dfrac{1}{2^{n}n!} {n \choose k} \dfrac{(\lambda^{2}-1)^{n-k}}{\lambda ^{n-2k}} \dfrac{d^{n}}{dx^{n}} (x^{2}-1)^{k}.
\label{Eq7}
\end{equation}
Recognizing the Rodrigues representation of the Legendre polynomial of degree $k$ in the previous equation and substituting for $P_{k}(x)$, we obtain:
\begin{equation}
P_{n}(\lambda x) = \sum_{k=\lceil n/2 \rceil}^{n} \dfrac{1}{2^{n-k}(n-k)!} \dfrac{(\lambda^{2}-1)^{n-k}}{\lambda ^{n-2k}} \dfrac{d^{n-k}}{dx^{n-k}} P_{k}(x),
\label{Eq8}
\end{equation}
which can be rewritten as:
\begin{equation}
P_{n}(\lambda x) = \sum_{k=0}^{\lfloor n/2 \rfloor}  a_{\lambda,n,k} \dfrac{d^{k}}{dx^{k}} P_{n-k}(x),
\label{Eq9}
\end{equation}
with
\begin{equation}
a_{\lambda,n,k}=\dfrac{\lambda ^{n-2k}(\lambda^{2}-1)^{k}}{2^{k}(k)!}.
\label{Eq9b}
\end{equation}
Equation (\ref{Eq9}) shows that $P_{n}(\lambda x)$ can be expressed as a sum of multiple derivatives of Legendre polynomials, which, in turn, can be expanded as a sum of Legendre Polynomials, as proved in the following section. It will be shown that:
\begin{equation}
\dfrac{d^{k}}{dx^{k}}P_{n-k}(x)=\sum_{i=0}^{\lfloor (n-2k)/2 \rfloor} \alpha_{n-2k-2i}P_{n-2k-2i}(x)
\label{Eq10}
\end{equation}
with
\begin{eqnarray}
\alpha_{n-2k-2i} &=&  \dfrac{2^{k+2i}(n-k-\nicefrac{1}{2})^{\underline{k}}(n-k-i)^{\underline{i}}(n-2k-\nicefrac{1}{2})^{\underline{2i}}}{(2i)^{\underline{2i}}(n-k-\nicefrac{1}{2})^{\underline{i}}} \nonumber \\
&& \qquad \qquad \qquad -\sum_{l=0}^{i-1} \dfrac{(2(n-2k-i-l))^{\underline{2(i-l)}}}{\left(2(i-l)\right)^{\underline{2(i-l)}}} \alpha_{n-2k-2l}.
\label{Eq11}
\end{eqnarray}
Combining equations (\ref{Eq9}), (\ref{Eq9b}), and (\ref{Eq10}), it follows that:
\begin{equation}
P_{n}(\lambda x) = \sum_{k=0}^{\lfloor n/2 \rfloor} \dfrac{\lambda ^{n-2k}(\lambda^{2}-1)^{k}}{2^{k}(k)!}  \sum_{i=0}^{\lfloor (n-2k)/2 \rfloor} \alpha_{n-2k-2i}P_{n-2k-2i}(x).
\label{Eq12}
\end{equation}
Finally, after rearranging the terms in the second sum of the previous equation, we obtain the following:
\begin{equation}
P_{n}(\lambda x) = \sum_{k=0}^{\lfloor n/2 \rfloor} b_{\lambda,n,k} P_{n-2k}(x) 
\label{Eq13}
\end{equation}
where,
\begin{equation}
b_{\lambda,n,k}= \sum_{i=0}^{\textrm{max}\{(k-1),0\}} \dfrac{\lambda ^{n-2k+2i}}{2^{k-i}(k-i)!}(\lambda^{2}-1)^{k-i}\alpha_{n,k,i},
\label{Eq14}
\end{equation}
and 
\begin{eqnarray}
\alpha_{n,k,i} &=&  \dfrac{2^{k+i}(n-k+i-\nicefrac{1}{2})^{\underline{k-i}}(n-k)^{\underline{i}}(n-2k+2i-\nicefrac{1}{2})^{\underline{2i}}}{(2i)^{\underline{2i}}(n-k+i-\nicefrac{1}{2})^{\underline{i}}} \nonumber \\
&& \qquad \qquad \qquad -\sum_{\substack{l=0 \\ i \neq 0}}^{i-1} \dfrac{(2(n-2k+i-l))^{\underline{2(i-l)}}}{\left(2(i-l)\right)^{\underline{2(i-l)}}} \alpha_{n,k-i+l,l}.
\label{Eq15}
\end{eqnarray}
 
 \subsection{Multiple derivative of Legendre polynomials}
 
From the recurrence relations \cite{Polyanin07}:
\begin{eqnarray}
\label{Eq16}
&& (n+1)P_{n+1}(x)-(2n+1)xP_{n}(x)+nP_{n-1}(x)=0, \qquad \textrm{and} \\
&& (x^{2}-1)\dfrac{d}{dx}P_{n}(x)=n \left[ xP_{n}(x) - P_{n-1}(x) \right] \nonumber,
\end{eqnarray}
it can be shown that:
\begin{equation}
\dfrac{d}{dx}\left[ P_{n+1}(x)-P_{n-1}(x)\right]=(2n+1)P_{n}(x).
\label{Eq17}
\end{equation}
From the previous equation, it follows that:
\begin{equation}
\dfrac{d}{dx}P_{n}(x)=\dfrac{2P_{n-1}(x)}{||P_{n-1}||^{2}}+\dfrac{2P_{n-3}(x)}{||P_{n-3}||^{2}}+\dfrac{2P_{n-5}(x)}{||P_{n-5}||^{2}}+...,
\label{Eq18}
\end{equation}
where $||P_{n}||=\sqrt{2/(2n+1)}$. By recursively deriving equation (\ref{Eq18}), we finally obtain:
\begin{equation}
\dfrac{d^{k}}{dx^{k}}P_{n}(x)=\sum_{i=0}^{\lfloor (n-k)/2 \rfloor} \alpha_{n-k-2i}P_{n-k-2i}(x)
\label{Eq19}
\end{equation}
The expansion coefficients $\alpha_{n-k-2i}$ can be determined by replacing the Legendre polynomials on each side of the previous equation with their corresponding hypergeometric series \cite{Gradshteyn14} and matching the coefficients of each term. Using Murphy's formula, Legendre polynomials can be written as:
\begin{equation}
P_{n}(x)= {}_{2}F_{1}(-n,n+1,1;(1-x)/2) = \sum_{j=0}^{\infty} \dfrac{(-n)_{j}(n+1)_{j}}{(1)_{j}}\dfrac{\left[(1-x)/2\right]^{j}}{j!},
\label{Eq20}
\end{equation}
where $(a)_{j}$ denotes the rising factorial. By setting $z=(1-x)/2$, and recognizing that $dz^{k}=(-1/2)^{k}dx^{k}$, the $k^{th}$ derivative of the Legendre polynomial of degree $n$ with respect to the variable $x$  is given by:
\begin{eqnarray}
\label{Eq21}
\dfrac{d^{k}}{dx^{k}}P_{n}(x) &=& \dfrac{d^{k}}{dz^{k}} \dfrac{dz^{k}}{dx^{k}} {}_{2}F_{1}(-n,n+1,1;(1-x)/2) \nonumber \\
                                              &=& \left(-\dfrac{1}{2}\right)^{k} \dfrac{d^{k}}{dz^{k}} {}_{2}F_{1}(-n,n+1,1;z) \nonumber \\
                                              &=& \left(-\dfrac{1}{2}\right)^{k}  \dfrac{(-n)_{k}(n+1)_{k}}{k!}{}_{2}F_{1}(-n+k,n+1+k,1+k;z) \nonumber\\
                                              &=& {n \choose k} \dfrac{(n+1)_{k}}{2^{k}}{}_{2}F_{1}(-n+k,n+1+k,1+k;z).
\end{eqnarray}
By recognizing that the rising factorial $(-n+k)_{i}$ is equal to zero for $i \geq n-k+1$, it then follows that:
\begin{equation}
\dfrac{d^{k}}{dx^{k}}P_{n}(x)= {n \choose k} \dfrac{(n+1)_{k}}{2^{k}} \sum_{j=0}^{n-k}\dfrac{(-n+k)_{j}(n+1+k)_{j}}{(1+k)_{j}}\dfrac{z^{j}}{j!}.
\label{Eq22}
\end{equation}
The same reasoning can be applied to the right side of equation (\ref{Eq19}), leading to:
\begin{eqnarray}
\label{Eq23}
P_{n-k-2i}(x) &=& {}_{2}F_{1}(-n+k+2i, n-k-2i+1,1;z) \nonumber \\
                    &=& \sum_{j=0}^{n-k-2i}\dfrac{(-n+k+2i)_{j}(n-k-2i+1)_{j}}{(1)_{j}}\dfrac{z^{j}}{j!} 
\end{eqnarray}
By matching the coefficients of each term in equations (\ref{Eq22}) and (\ref{Eq23}), we obtain a set of $(n-k)$ coupled equations:
\begin{eqnarray}
\label{Eq24}
A_{n-k} &=& \alpha_{n-k}B_{n-k, n-k} \\
A_{n-k-1} &=& \alpha_{n-k}B_{n-k-1, n-k} \nonumber \\
A_{n-k-2} &=& \alpha_{n-k}B_{n-k-2, n-k} + \alpha_{n-k-2}B_{n-k-2, n-k-2} \nonumber \\
A_{n-k-3} &=& \alpha_{n-k}B_{n-k-3, n-k} + \alpha_{n-k-2}B_{n-k-3, n-k-2} \nonumber \\
A_{n-k-4} &=& \alpha_{n-k}B_{n-k-4, n-k} + \alpha_{n-k-2}B_{n-k-4, n-k-2} + \alpha_{n-k-4}B_{n-k-4, n-k-4} \nonumber \\
A_{n-k-5} &=& \alpha_{n-k}B_{n-k-5, n-k} + \alpha_{n-k-2}B_{n-k-5, n-k-2} + \alpha_{n-k-4}B_{n-k-5, n-k-4} \nonumber \\
.... \nonumber
\end{eqnarray}
where, for the sake of clarity, we have set:
\begin{eqnarray}
\label{Eq25}
A_{j} &=& {n \choose k} \dfrac{(n+1)_{k}}{2^{k}}\dfrac{(-n+k)_{j}(n+1+k)_{j}}{(1+k)_{j}} \\
B_{j,n-k} &=& \dfrac{(-n+k)_{j}(n-k+1)_{j}}{(1)_{j}} \nonumber \\
B_{j, n-k-2} &=& \dfrac{(-n+k+2)_{j}(n-k-1)_{j}}{(1)_{j}} = B_{j,n-k}\dfrac{(n-k-j)^{\underline{2}}}{(n-k+j)^{\underline{2}}} \nonumber \\
B_{j, n-k-4} &=& \dfrac{(-n+k+4)_{j}(n-k-3)_{j}}{(1)_{j}} = B_{j,n-k}\dfrac{(n-k-j)^{\underline{4}}}{(n-k+j)^{\underline{4}}} \nonumber \\
.... \nonumber \\
B_{j, n-k-2i} &=& \dfrac{(-n+k+2i)_{j}(n-k-2i+1)_{j}}{(1)_{j}} = B_{j,n-k}\dfrac{(n-k-j)^{\underline{2i}}}{(n-k+j)^{\underline{2i}}} \nonumber \\
.... \nonumber 
\end{eqnarray}
In the previous equation, the symbol $(x)^{\underline{n}}$ is used to represent the falling factorial. From equations (\ref{Eq24}) and identities (\ref{Eq25}), a recurrence relation for the coefficients $\alpha_{n-k-2i}$ can be derived. It can be shown that:
\begin{equation}
\alpha_{n-k-2i} =  \dfrac{2^{k+2i}(n-\nicefrac{1}{2})^{\underline{k}}(n-i)^{\underline{i}}(n-k-\nicefrac{1}{2})^{\underline{2i}}}{(2i)^{\underline{2i}}(n-\nicefrac{1}{2})^{\underline{i}}}-\sum_{l=0}^{i-1} \dfrac{(2(n-k-i-l))^{\underline{2(i-l)}}}{\left(2(i-l)\right)^{\underline{2(i-l)}}} \alpha_{n-k-2l}.
\label{Eq26}
\end{equation}

\section{Conclusion}

In this work, we have studied the scaling properties of Legendre polynomials $P_{n}(x)$. We have demonstrated that $P_{n}(\lambda x)$, where $\lambda$ is a constant, can be expanded either as a sum of Legendre polynomials $P_{n}(x)$ or their multiple derivatives, and we have obtained an explicit expression for the expansion coefficients. In addition, we have shown that the multiple derivative $d^{k}P_{n}(x)/dx^{k}$ can also be expressed as a sum of Legendre polynomials and we derived a recurrence relation for the coefficients.

\section*{Acknowledgements}
This work was supported by the U.S. Department of Energy, Office of Science, Basic Energy Sciences, under Award No. DE-SC0017984.

\end{document}